\documentclass[final,11pt,leqno]{amsart}

\usepackage{amsmath,amssymb,amscd,amsthm}
\usepackage{latexsym}
\usepackage{graphicx}

\usepackage{showkeys}

\newtheorem{theorem}{Theorem}
\newtheorem{lemma}{Lemma}
\newtheorem{proposition}{Proposition}

\newtheorem{definition}{Definition}
\newtheorem{corollary}{Corollary}

\newtheorem{claim}{Claim}

\newcommand{\f}[2]{\frac{#1}{#2}}
\newcommand{\dpr}[2]{\langle #1,#2 \rangle}

\newcommand{\de}{\delta}

\newcommand{\ka}{\kappa}
\newcommand{\la}{\lambda}

\newcommand{\si}{\sigma}

\newcommand{\vp}{\varphi}
\newcommand{\om}{\omega}

\newcommand{\rone}{\mathbf R^1}

\newcommand{\ca}{\mathcal A}
\newcommand{\cL}{\mathcal L}
\newcommand{\cJ}{\mathcal J}

\newcommand{\p}{\partial}

\newcommand{\beq}{\begin{equation}}
\newcommand{\eeq}{\end{equation}}
\newcommand{\beqna}{\begin{eqnarray*}}
\newcommand{\eeqna}{\end{eqnarray*}}
\newcommand{\beqn}{\begin{equation*}}
\newcommand{\eeqn}{\end{equation*}}
\newcommand{\bp}{\begin{proof}}
\newcommand{\ep}{\end{proof}}
\newcommand{\bprop}{\begin{proposition}}
\newcommand{\eprop}{\end{proposition}}
\newcommand{\bt}{\begin{theorem}}
\newcommand{\et}{\end{theorem}}
\newcommand{\bex}{\begin{Example}}
\newcommand{\eex}{\end{Example}}
\newcommand{\bc}{\begin{corollary}}
\newcommand{\ec}{\end{corollary}}
\newcommand{\bcl}{\begin{claim}}
\newcommand{\ecl}{\end{claim}}
\newcommand{\bl}{\begin{lemma}}
\newcommand{\el}{\end{lemma}}

\begin{document}

\title
[ Transverse instability for periodic waves] {Transverse
instability for periodic waves of KP-I and Schr\"odinger
equations}

\author{Sevdzhan Hakkaev}
\author{Milena Stanislavova}
\author{Atanas Stefanov}

\address{Sevdzhan Hakkaev
Faculty of Mathematics and Informatics, Shumen University, 9712
Shumen, Bulgaria}\email{shakkaev@fmi.shu-bg.net}
\address{Milena Stanislavova
Department of Mathematics, University of Kansas, 1460 Jayhawk
Boulevard,  Lawrence KS 66045--7523} \email{stanis@math.ku.edu}
\address{Atanas Stefanov
Department of Mathematics, University of Kansas, 1460 Jayhawk
Boulevard,  Lawrence KS 66045--7523}

\email{stefanov@math.ku.edu}

\thanks{
Stanislavova  supported in part by NSF-DMS \# 0807894. Stefanov
supported in part by NSF-DMS \# 0908802 .}
\date{\today}

\subjclass[2000]{Primary }

\keywords{transverse instability, periodic traveling waves, KP-I equation, Schr\"odinger equation}

\begin{abstract}
We consider the  quadratic and cubic
KP  - I  and  NLS models in $1+2$ dimensions with periodic boundary conditions.
 We show that the spatially periodic travelling waves (with period $K$) in the form $u(t,x,y)=\vp(x-c t)$ are spectrally and linearly  unstable, when the perturbations are taken to be with the same period. This strong instability implies other instabilities considered recently - for example  with respect to perturbations with  periods $nK, n=2, 3, \ldots$ or bounded perturbations.
\end{abstract}

\maketitle

\section{Introduction and Statements of main results}

The existence and stability properties of special solutions of nonlinear differential equations is   an important question both from theoretical and practical point of view.
Many equations  describing
 wave motion typically  feature  traveling wave solutions.   The problem of the
orbital stability of solitary waves for nonlinear dispersive
equations goes back to the works of Benjamin \cite{Be1} and Bona
\cite{Bo}. 
Another approach is to linearize the equation around the solitary wave
and look for linear stability based on the spectrum of the linear solution operator.
Extending the ODE ideas to partial differential equations has introduced
 a number of new issues. In infinite dimensions, the relation between the linearization and the full
nonlinear equations is far more complicated.
Another nontrivial issue arises at the linear level, since all of
the known proofs for the existence of invariant manifolds are based upon
the use of the solution group (or semigroup) generated by the linearization.
However, in any actual problem, the information available will, at best, be of
the spectrum of the infinitesimal generator, that is, the linearized
equation and not its solution operator. Relating the spectrum of the
infinitesimal generator to that of the group is a spectral mapping problem
that is often non-trivial. All of these three problems - 
 spectral stability, linear stability and nonlinear stability, 
 have been extensively studied for solitary wave solutions.

While the existence and stability of such solutions on  the whole
space case has been well-studied, the questions about existence and stability of spatially periodic traveling waves have not received much attention until recently.
One of the first results on stability of periodic solutions of the Korteweg-de Vries(KdV)
equation was obtained by McKean \cite{M}. Based on the
integrability of the KdV equation the stability of all periodic
finite-genus solutions has been established. Recently Angulo, Bona and
Scialom \cite{ABS} investigated the orbital stability of cnoidal waves for the
KdV equation with respect to perturbations of the same period.
The linear stability/instability of some of these solutions with respect to different types of perturbations
has been developed in the last couple of years, see for example  \cite{HK}, \cite{BJ} and \cite{BD}.
Other new explicit formulae for periodic traveling waves of dnoidal
type   together with their stability have been obtained in
\cite{An2, AnNa, HIK1, HIK2, HIK3}.

An  interesting aspect of the theory is when one considers
the one-dimensional waves as solutions in  two-dimensional models.
One generally refers to this as the question for {\it transverse stability} of such waves.
 The transverse stability of traveling waves  is
associated with a class of perturbations   traveling
 transversely to the direction of the basic traveling wave.

The problem of transverse stability/instability of solitary waves goes back to a work 
 by Kadomtsev and Petviashvili \cite{KP} for
 KdV solitary waves.  It turns out that solitary waves are transverse
stable in the case of KP-I and transverse instable in the case of
KP-II.

Recently, Rousset and Tzvetkov, \cite{RT1, RT2} provided general
criterium for transverse instability for traveling waves of
Hamiltonian partial differential equations, which was then applied
to various examples. Johnson and Zumbrun \cite{JZ} investigated
the stability  {\it periodic}  traveling waves of the generalized
KdV equation to two dimensional perturbations, which are
nonperiodic (bounded) in the generalized KP equation and have long
wavelength in transverse direction. By analyzing high and low
frequency limits of the appropriate periodic Evans function they
derived an instability criterion for the transverse instability.
This criteria  is then applied to the KdV and modified KdV
equations. The authors proved that the periodic traveling waves of
the KdV equation are unstable to long wavelength transverse
perturbations and that cnoidal, and dnoidal traveling waves for
modified KdV equation are transverse unstable to long wavelength
perturbations in KP-II and KP-I respectively. Haragus \cite{H}
considered the transverse spectral stability of {\it small
periodic} traveling wave solutions of the KdV equation with
respect to perturbations in KP-I and KP-II which are either
periodic in the direction of perturbation or nonperiodic
(localized or bounded) and have long wavelength in the transverse
direction.

 In this paper, we  prove transverse instability of certain periodic solutions  of the
 Kadomtsev-Petviashvili-I equation and the nonlinear Schr\"odinger equation.
 More precisely, we consider  periodic traveling waves of
the KdV and mKdV equation, which in turn also solve the KP-I equation, 
while our second example concerns spatially periodic standing waves of the non-linear Schr\"odinger equation (NLS).
Before we continue with the specifics of our results, we 
 outline the general scheme and we  give some definitions.

In this paper we only deal with the stability information provided by the linearized equation\footnote{i.e. we will not consider the full non-linear equation satisfied by $v$, which would of course amount to non-linear stability/instability results.}. Suppose that the linearized equation is in the form of an evolution equation
\begin{equation}
\label{evol}
v_t=\ca v.
\end{equation}
We use the following definition of spectral and linear  stability
\begin{definition}
\label{defi:1}
Assume that $\ca=\ca(\vp)$  generates a $C_0$ semigroup on a Banach space $X$.
We say that the solution $\vp$ with linearized problem \eqref{evol} is spectrally stable, if
$\si(\ca)\subset \{\la: \Re\la\leq 0\}$.

We say that the the solution $\vp$ with linearized problem \eqref{evol} is linearly
 stable, if the growth bound for the semigroup $e^{t\ca}$ is non-positive.
 Equivalently, we require that every solution of \eqref{evol} with $v(0)\in X$ has the property
 $$
 \lim_{t\to \infty} e^{-\de t} \|v(t, \cdot)\|=0
 $$
 for every $\de>0$.
\end{definition}
{\bf Remarks:}
 We recall that by the spectral mapping theorem for point spectrum
 $\si_{p.p.}(e^{t\ca})\setminus \{0\}
= e^{t\si_{p.p.}(\ca)}$. There is however only the inclusion
$\si_{ess}(e^{t\ca})\setminus \{0\} \supseteq e^{t\si_{ess}(\ca)}$, which is the reason that one canoot, in general (and in the absence of the so-called spectral mapping theorem), deduce linear stability from spectral stability.
In fact, due to the spectral inclusions above, linear stability implies spectral stability, but in general the converse is false.

However, in the cases considered in this paper the spectrum consists entirely of eigenvalues and the two notions are equivalent (since there is a spectral mapping theorem for eigenvalues, as indicated above).  Thus, we will concentrate on the spectral stability from now on.

 \subsection{KP - I equation}   Consider
 the spatially periodic  KP - I equation
  \begin{equation}
  \label{kp}
     \left\{\begin{array}{l}
     (u_t+\partial_{xxx}u+\p_x(f(u)))_x- \partial_{yy}u=0, \ \
     (t,x,y)\in \rone_+\times [0,K_1]\times [0,K_2] \\
     u(t,x+K_1,y)=u(t,x,y); u(t,x,y+K_2)=u(t,x,y)
     \end{array}
     \right.
   \end{equation}
   where $f$ is smooth function\footnote{We only consider the cases
   $f(u)=u^2, \pm u^3$, but other choices certainly make sense mathematically.}.  It is  known that solutions exists, at least locally, when the data is in the product Sobolev spaces  $f\in H^{3,3}([0,K_1]\times [0,K_2])$,   see for example \cite{IMS}.

  In this paper, we will be interested in the stability properties of a class of special solutions, namely the  periodic  traveling waves solution of the modified KdV equation. That is, we look for solutions  in the form $v(t,x)=\vp(x-ct)$, $\vp(x+K_1)=\vp(x)$,  so that
  $$
  v_t+\partial_{xxx}v+\p_x(f(v))=0, \ \  x\in [0,K_1].
  $$
 Clearly then $u(t,x,y):=\vp(x-ct)$ is a solution of the KP - I equation
 \eqref{kp}.  We construct these solutions $\vp_c$  explicitly in Section \ref{sec:2} below.
Periodic travelling-wave solution  are determined from Newton's
equation which we will write below in the form
$\varphi'^2=U(\varphi)$. Therefore by using the well-known
properties of the phase portrait of Newton's equation in the
$(\varphi,\varphi')$-plane, one can establish that under fairly
general conditions,  that there exists a family of periodic
solutions $\varphi(y)=\varphi(c,\varphi_0;y)$  and $\varphi_0
=\min \varphi$.  Moreover, if $T=T(c,\varphi_0)$ (in particular,
their period turns out to depend on the speed parameter $c$ and
an elliptic modulus $\kappa$) is the minimal (sometimes called {\it
fundamental}) period of $\varphi$, then $\varphi$ has exactly one
local minimum and one local maximum in $[0,T)$. Therefore
$\varphi'$ has just two zeroes in each semi-open interval of
length $T$. By Floquet theory, this means that $\varphi'$ is
either the second or the third eigenfunction of the periodic
eigenvalue problem.

 In order to explain the stability/instability results, we need to  linearize the equation \eqref{kp} about the periodic traveling wave solution. Namely, write an ansatz in the form
  $u(t,x,y)=\varphi(x-ct)+v(t, x-ct,y)$, which we plug in \eqref{kp}.
  After ignoring all nonlinear in $v$ terms, we arrive at the following linear equation for $v$
     \begin{equation}\label{le}
       (v_t  + v_{xxx}-c v +(f'(\varphi)v)_x)_x-\partial_{yy}v=0.
     \end{equation}
 If the variable  $v$ has the mean-zero property in $x$ (i.e.  $\int_0^{K_1} v(t,x,y) dx=0$), then one may invert the  operator $\p_x$ (by defining $(\p_x^{-1} f)(x): =\int_0^x f(y)dy$) and thus recast \eqref{le} in the evolution equation form
 \begin{equation}\label{lemod}
       v_t= \p_x( -\p_x^2 +c  -f'(\varphi)) v+\p_x^{-1} \partial_{yy}v
     \end{equation}
The question for stability/instability of traveling wave  solutions of the KP - I equation has attracted a lot of attention in the last few years (see \cite{JZ}, \cite{H}).

As we have indicated above, we  restrict our attention to
spectral considerations for the generator. In order to establish instability, we  seek
 solutions in the form
   $$
   v(t,x,y)=e^{\sigma t}e^{iky}V(x),
   $$
 where $\sigma \in \mathbb{C}, \; k\in \mathbb{R}$ and $V(x)$ is
   periodic function with same period as the periodic traveling wave
   solution $\varphi(x)$. Clearly, such solutions will be also periodic in the $y$ variable, with period $K_2=2\pi/k$. Thus, if we   manage to show existence of such $V=V(\si, k)$ with some
   $\si>0$, we will have shown transeverse spectral instability of the traveling wave
   solution $\varphi(x)$.

  We further specialize  $V$ in the form
   $V=\partial_x U$. Plugging  in \eqref{lemod}  yields the equation
      $$
      -\sigma \partial_x U=(-\partial_x
      (-\partial_{xx}+c-f'(\varphi))\partial_x+k^2)U.
      $$
 This eigenvalue problem is therefore in the form
   \begin{equation}\label{2.1}
    \sigma A(k)U=L(k)U,
  \end{equation}
  with
     $$
     A(k)=-\partial_x, \; \;
     L(k)=-\partial_x(-\partial_{xx}+c-f'(\varphi))\partial_x+k^2.
     $$
where $L(k), A(k)$ are operators which depend on the real
parameter $k$ on some Hilbert space  $H$.

\subsection{The Nonlinear Schr\"odinger Equation}

Another object of investigation will be the spacially periodic solutions of the Nonlinear Schr\"odinger Equation (NLS).
\begin{equation}
\label{sch}
\left\{\begin{array}{l}
i u_t-(u_{xx}+u_{yy}) - f(|u|^2) u=0, (t,x,y)\in \rone_+\times [0,K_1]\times [0,K_2]\\
u(t,x+K_1, y)=u(t,x, y);\  u(t,x, y+K_2)=u(t,x, y).
\end{array}
\right.
\end{equation}
where $f$ is a smooth function.
Looking for standing waves in the form $u(t,x)=e^{-i \om t}\vp(x)$ results in the ordinary differential equation
\begin{equation}
\label{m1}
\om \vp-\vp''-f(\vp^2)\vp=0.
\end{equation}
We now derive the linearized equation for small perturbation of the wave $e^{-i \om t} \vp$. Write the ansatz $u=e^{-i \om t} (\vp+v(t,x,y))$. For the nonlinear term, we have
\begin{eqnarray*}
& & f(|u|^2)= f(|\vp+v|^2)=f(\vp^2+2\vp\Re v+|v|^2)=f(\vp^2)+2 f'(\vp^2) \vp\Re v+O(v^2).
\end{eqnarray*}
We get, after disregarding $O(v^2)$ terms and taking into account \eqref{m1},
$$
i v_t +\om v-(v_{xx}+v_{yy})-f(\vp^2)v -2 f'(\vp^2) \vp^2 \Re v=0.
$$
We are looking for unstable solutions in the form $v(t,x,y)=e^{\si t} \cos(k y) V(x)$, where $V$ is a complex-valued function. We obtain
$$
i \si V+\om V-V''+k^2 V -f(\vp^2) V-2f'(\vp^2)\vp^2 \Re V=0
$$
Let $V=v_1+i v_2$, where $v_1, v_2$ are real-valued functions.  This  gives the following system for $v_1, v_2$
\begin{eqnarray*}
& & \si v_1-v_2''+\om v_2+k^2 v_2-f(\vp^2)v_2=0\\
& &-\si v_2-v_1''+\om v_1+k^2 v_1-f(\vp^2) v_1-2 f'(\vp^2) \vp^2 v_1=0.
\end{eqnarray*}
Denote
\begin{eqnarray*}
& &   \cL_+=-\p_x^2+\om-f(\vp^2), \\
& &   \cL_-=-\p_x^2+\om -f(\vp^2)-2f'(\vp^2)\vp^2.
\end{eqnarray*}
This allows us to write the linearized problem as follows
\begin{equation}
\label{m2}
\si\left(\begin{array}{c}
v_1 \\ v_2
\end{array}\right)+ \left(\begin{array}{c c}
0 & \cL_++k^2\\  -(\cL_-+k^2) & 0
\end{array}\right)\left(\begin{array}{c}
v_1 \\ v_2
\end{array}\right)=0.
\end{equation}
Let  $J=\left(\begin{array}{c c}
0 & 1 \\ -1 & 0
\end{array}\right)$ and $\left(\begin{array}{c}
v_1 \\ v_2
\end{array}\right)=J \left(\begin{array}{c}
z_1 \\ z_2
\end{array}\right)$. Note that $J^*=J^{-1}=-J$. In terms of $z_1, z_2$, we have the equation
$$
\si J \left(\begin{array}{c}
z_1 \\ z_2
\end{array}\right)=-J \left(\begin{array}{c c}
\cL_-+k^2 & 0\\  0 & \cL_++k^2
\end{array}\right) J \left(\begin{array}{c}
z_1 \\ z_2
\end{array}\right)
$$
Thus, we have managed to recast the problem in the form \eqref{2.1}, this time with
\begin{equation}
\label{m3}
A(k)=\si J; \ \ L(k)= -J \left(\begin{array}{c c}
\cL_-+k^2 & 0\\  0 & \cL_++k^2
\end{array}\right) J = -J \left(\begin{array}{c c}
\cL_-& 0\\  0 & \cL_+
\end{array}\right) J+k^2 Id
\end{equation}
{\bf Remark:} We would like to give  the important case $f(z)=\sqrt{z}$ some more consideration, due to the fact that the function $\sqrt{z}$ fails to be differentiable at zero. Nevertheless, we still have
$$
\sqrt{|\vp+v|}=\vp+\Re v+o(v),
$$
and we still obtain the formula
\begin{eqnarray*}
& &   \cL_+=-\p_x^2+\om-\vp, \\
& &   \cL_-=-\p_x^2+\om -2\vp,
\end{eqnarray*}
as we would, if we were to use the derivative of the function $f(z)$ in the generic definition of $\cL_\pm$ above. The difference of course is in the fact that the remainder term is only $o(v)$ instead of $O(v^2)$, but this of course is irrelevant for the linear theory that we develop here.
\subsection{Main Results}
Our first result concerns the transverse instability of the cnoidal solutions of the
KP - I equation.
\begin{theorem}(transverse instability for cnoidal solutions of KP - I)\\
\label{theo:2} Considering the KP - I equation (i.e. \eqref{kp}
with $f(u)=\f{u^2}{2}$).  It supports the cnoidal solutions given
by \eqref{3.13} below. Then, there exists a period $K_2$ depending on
the particular cnoidal solution, so that the cnoidal waves are spectrally and
linearly unstable for all values of the parameters $\kappa \in
(0,1)$ and $T$ given by \eqref{a15}.
\end{theorem}

Next, we state our main result regarding  transverse instability of the dnoidal solutions of the
modified KP - I equation.
\begin{theorem}(transverse instability for dnoidal solutions of modified KP - I)\\
\label{theo:5} Consider the modified KP - I equation, that is
\eqref{kp} with $f(u)=u^3$. Then, there exists a period $K_2$ depending on
the particular dnoidal solution, so that the  dnoidal solutions described by
\eqref{3.24} below are spectrally and linearly unstable for all
values of the parameters $\kappa \in (0,1)$ and the corresponding $T$.
\end{theorem}
Finally, we have the following result, which shows transverse instability for standing waves of the quadratic and cubic NLS. That is, we shall be considering
\eqref{sch} with $f(z)=\sqrt{z}$ and $f(z)=z$.
\begin{theorem}(transverse instability for standing wave solutions of  NLS)\\
\label{theo:10}
The quadratic (focussing) Schr\"odinger equation\footnote{i.e. $f(z)=\sqrt{z}$}  \eqref{sch} admits cnoidal solutions in the form \eqref{s1.5}. There exists $K_2$, depending on the specific solution, so that these solutions are spectrally and linearly unstable for all values of the parameter $\ka\in (0,1)$.

The cubic (focussing) Schr\"odinger equation\footnote{i.e. $f(z)=z$}  \eqref{sch} supports dnoidal solutions in the form \eqref{s2.4}. There exists $K_2$, depending on the specific solution, so that these solutions are spectrally and linearly unstable for all values of the parameter $\ka\in (0,1)$.
\end{theorem}
{\bf Remarks:}
\begin{itemize}
\item As a consequence of the three theorems above, one may deduce spectral instability, when the perturbations are taken to be periodic (with period equal to integer times the period of the wave) or bounded functions.
\item Our method for showing transverse  instability fails for periodic snoidal waves of the defocussing modified KP - I equation, see Chapter \ref{ch:5}.  Beyond the technical issues, which prevents the relevant inequality \eqref{a:10} from being satisfied,  it would be interesting to further investigate the transverse stability/instability of these interesting waves.
\end{itemize}

\subsection{General instability criteria}
In our proofs, we use the following sufficient condition for instability.
\begin{theorem}
\label{theo:1} Assume that the operator $L(k)$
satisfies\footnote{Hereafter, we use the notation
$L'(k):=\f{d}{dk} L(k)$.}
\begin{enumerate}
  \item there exists $k_0>0$, so that
  $dim\  Ker[L(k_0)]=1$, say $Ker[L(k_0)]=span\{\vp\}$.
  \item $L'(k_0)\vp\neq 0$.
\end{enumerate}
Then, the equation \eqref{2.1} has a solution $U$ for some $k$,
sufficiently close to $k_0$ and for some sufficiently small
$\si>0$.
 In fact, there exists a continuous scalar
 function $k(\si): k(0)=k_0$ and a continuous $H$-valued function
 $U(\si): U(0)=\vp$, so that
 $$
 \si A(k(\si)) U(\si)=L(k(\si)) U(\si),
 $$
 for all $0<\si<<1$.
\end{theorem}
{\bf Note:} This is a variant of a theorem used by Groves-Haragus and Sun, \cite{GHS}. The interested reader should also explore the simple exposition in \cite{RT},  where several examples about transverse instability on the whole space are worked out in detail using the same techniques.
\begin{proof}
We quickly indicate the main ideas of the proof. \\
  Let  $U=\varphi
  +V$, with
    $$V\in \varphi^{\perp}=\{ V\in H, \; \; (V, \varphi)=0 \}.$$

    Consider the equation $G(V, k, \sigma)=0$, with $\sigma>0$ and
      $$G(V, k, \sigma)=L(k)\varphi +L(k)V-\sigma A(k)\varphi
      -\sigma A(k)V.$$
    We have
      $$
      \dpr{D_{V,k}(0, k_0, 0)}{[\omega, \mu]}=\mu L'(k_0)\varphi+L(k_0)\omega$$
      and
      $D_{V,k}(0, k_0, 0)$ is a bijection from
      $\varphi^{\perp}\times \mathbb{R}$ to $H$. Thus from the
      implicit function theorem follows that for $\sigma$ in a
      neighborhood of zero there exists $k(\sigma)$ and
      $V(\sigma)$ such that $G(V(\sigma), k(\sigma), \sigma)=0$.
   \end{proof}
Clearly, in view of Theorem \ref{theo:1} and the spectral problem
\eqref{2.1}, we will have proved Theorem \ref{theo:2} and Theorem \ref{theo:5}, provided we can verify the conditions $(1), (2)$ of Theorem \ref{theo:1} for the operator
$$
L(k)=-\p_x\cL\p_x+k^2= -\p_x(-\p_{xx}+c-f'(\vp))\p_x +k^2.
$$
Similarly, for Theorem \ref{theo:10}, due to the representation \eqref{m3}, it suffices  to
verify conditions $(1), (2)$ of Theorem \ref{theo:1} for the operator
$$
L(k)=\cJ^{-1} \cL \cJ+k^2=\cJ^{-1}\left(\begin{array}{cc} \cL_- &  0 \\
0 & \cL_+  \end{array} \right)\cJ+k^2.
$$
This clearly necessitates a somewhat detailed study of the spectral picture for  the operators $\cL, \cL_\pm$. Luckily, after one constructs the traveling/standing waves for our models in terms of elliptic functions, we will be able to obtain some information about the spectra of $\cL$ and $\cL_\pm$, which will allow us to check condition $(1)$ in Theorem \ref{theo:1}.

 The paper is organized as follows. In Section \ref{sec:2}, we construct the eigenfunctions. In Section \ref{sec:30}, we describe the structure of the first few eigenvalues, together with the associated eigenfunctions for $\cL$ and $\cL_\pm$. In section \ref{sec:44}, we give the proof of Theorem \ref{theo:2} by verifying conditions $(1), (2)$. This requires some spectral theory, together with the specific spectral information for $\cL, \cL\pm$, obtained in Section \ref{sec:2}.
 In section \ref{ch:5}, we show that an identical approach for
 the defocusing modifed KP - I equation  {\it fails} to give transverse instability. Thus, an interesting question is left open, namely - Are the snoidal solutions to theis problem transverse unstable?

  \section{Construction of periodic traveling waves}
  \label{sec:2}
   We are looking for a traveling-wave solution for the equation
      \begin{equation}\label{3.1}
        u_t+(f(u))_x+u_{xxx}=0
      \end{equation}
         of the form
$u(x,t)=\phi(x-ct)$. We assume that $\phi$ is smooth and bounded
in $\mathbb{R}$. The following two cases appear:

\vspace{1ex} (i) $\phi'\neq 0$ in $\mathbb{R}$ and
$\phi_-<\phi<\phi_+$ (corresponding to kink-wave solution);

\vspace{1ex} (ii) $\phi'(\xi)=0$  for some $\xi\in \mathbb{R}$.
Denote $\phi_0=\phi(\xi)$, $\phi_2=\phi''(\xi)$.

\vspace{2ex} \noindent Below we will deal with the second case.
Replacing in (\ref{3.1}) we get

\begin{equation}\label{3.2}
-c\phi'+(f(\phi))'+\phi'''=0.
\end{equation}
 Integrating (\ref{3.1}) twice, one obtains
 \begin{equation}\label{3.3}
 -c\phi+f(\phi)+\phi''=a
 \end{equation}

\begin{equation}\label{3.4}
\frac{\varphi'^2}{2}= b+a\phi+\frac{c}{2}\phi^2-F(\phi),\qquad
F(\phi)=\int_0^\phi f(s)ds
\end{equation}
 with some
constants $a, b$. In case (ii), one has respectively
$$\begin{array}{l}
a=f(\phi_0)-c\phi_0+\phi_2,\\
b=F(\phi_0)-\frac12c\phi_0^2-a\phi_0
=F(\phi_0)-\frac12c\phi_0^2-(-c\phi_0+f(\phi_0)+\phi_2)\phi_0.\end{array}
$$

Next we are going to look for periodic travelling-wave solutions
$\phi$. Consider in the plane $(X,Y)=(\phi,\phi')$ the Hamiltonian
system

\begin{equation}\label{3.5}
\begin{array}{l}
\dot{X}=Y=H_Y,\\
\dot{Y}=-f(X)+cX+a=-H_X,\end{array}
\end{equation}
 with a
Hamiltonian function
$$H(X,Y)=\frac{Y^2}{2}+F(X)-\frac{v}{2}X^2-aX.$$
Then (\ref{3.4}) becomes $H(\phi,\phi')=b$ and the curve $s\to
(\phi(s-s_0),\phi'(s-s_0))$ determined by (\ref{3.4}) lies on the
energy level $H=b$ of the Hamiltonian $H(X,Y)$. Within the
analytical class, system (\ref{3.5}) has periodic solutions if and
only if it has a center. Each center is surrounded by a continuous
band of periodic trajectories (called {\it period annulus}) which
terminates at a certain separatrix contour on the Poincar\'e
sphere. The critical points of center type of (\ref{3.5}) are
given by the critical points on $Y=0$ having a negative Hessian.
These are the points $(X_0,0)$ where:
\begin{equation}\label{3.6}
a+cX_0-f(X_0)=0,\;\;\; c-\phi'(X_0)<0.
\end{equation}
   (For simplicity, we will not consider here the case of a
degenerate center when the Hessian becomes zero.)

The above considerations lead us to the following statement.

\begin{proposition}
\label{prop:11}
Let $a$ and $c$
be constants such that conditions $(\ref{3.6})$ are satisfied for
some $X_0\in \mathbb{R}$. Then there is an open interval $\Delta$
containing $X_0$ such that:
\begin{enumerate}
\item[(i)]  For any $\phi_0\in \Delta$, $\phi_0<X_0$,
the solution of $(\ref{3.1})$ satisfying
$$\phi(\xi)=\phi_0,\quad  \phi'(\xi)=0,\quad
\phi''(\xi)=a+c\phi_0-f(\phi_0),$$ is periodic.
\item[(ii)]  If $\phi_1\in\Delta$, $\phi_1>X_0$ is the
nearest to $X_0$ solution of $H(X,0)=H(\phi_0,0)$, then
$\phi_0\leq\phi\leq\phi_1$.\\
\item[(iii)]  If $T$ is the minimal period of $\phi$,
then in each interval $[s,s+T)$,  the function $\phi$ has just one
minimum and one maximum ($\phi_0$ and $\phi_1$,
respectively) and it is strictly monotone elsewhere.
\end{enumerate}
\end{proposition}

\vspace{2ex} \noindent Denote
$$ U(s)=2b+2as+cs^2-2F(s)=2F(\phi_0)-c\phi_0^2-2a\phi_0+2as+cs^2-2F(s).$$ Then for
$\phi_0\leq\phi\leq \phi_1$ one can rewrite (\ref{3.4}) as
$\phi'(\sigma)=\sqrt{U(\varphi(\sigma))}$.  Integrating the
equation along the interval $[\xi,s]\subset [\xi,\xi+T/2]$ yields
an implicit formula for the value of $\phi(s)$:
\begin{equation}\label{3.7}
\int_{\phi_0}^{\phi(s)}\frac{d\sigma}{\sqrt{U(\sigma)}}=s-\xi,
\quad  s\in [\xi,\xi+T/2].
\end{equation}
For $s\in [\xi+T/2,\xi+T]$ one has $\varphi(s)=\varphi(T+2\xi-s)$.
We recall that the period function $T$ of a Hamiltonian flow
generated by $H_0\equiv \frac12Y^2-\frac12U(X)=0$ is determined
from
\begin{equation}\label{3.8}
T=\int_0^Tdt=\oint_{H_0=0}\frac{dX}{Y}=
2\int_{\varphi_0}^{\varphi_1}\frac{dX}{\sqrt{U(X)}}.
\end{equation}
 This is in fact the derivative (with respect to the
energy level) of the area surrounded by the periodic trajectory
through the point $(\phi_0,0)$ in the $(X,Y)=(\phi, \phi')$-plane.

\vspace{1ex} \noindent Consider the continuous family of
periodical travelling-wave solutions $\{u=\phi(x-ct)\}$ of
(\ref{3.1}) and (\ref{3.3}) going through the points
$(\phi,\phi')=(\phi_0,0)$ where $\phi_0\in\Delta^-$. For any
$\phi_0\in\Delta^-$, denote by $T=T(\phi_0)$ the corresponding
period. One can see (e.g. by using formula (\ref{3.7}) above) that
the period function $\phi_0\to T(\phi_0)$ is smooth. To check
this, it suffices to perform a change of the variable
\begin{equation}\label{3.9}
X=\frac{\phi_1-\phi_0}{2}s+\frac{\phi_1+\phi_0}{2}
 \end{equation}
in the integral (\ref{3.7}) and use that
\begin{equation}\label{3.10}
U(\varphi_0)=U(\varphi_1)=0.
\end{equation}
  Conversely, taking $c$,
$a$ to satisfy the conditions of Proposition \ref{prop:11} and fixing $T$ in a
proper interval, one can determine $\phi_0$ and $\phi_1$ as smooth
functions of $c$, $a$ so that the periodic solution $\phi$ given
by (\ref{3.7}) will have a period $T$. The condition for this is
the monotonicity of the period (for more details see).

\section{Spectral properties of the operators $\cL$ and $\cL_{\pm}$}
\label{sec:30}
We first construct the spectral representation of the KdV equation
\subsection{The operator $\cL$ for KdV}
\label{sec:31}
Consider the Korteweg-de Vries equation
\begin{equation}\label{KdV}
u_t+uu_x+u_{xxx}=0, \end{equation}
 which is a particular case of
(\ref{3.1}) with $f(u)=\frac{u^2}{2}$. In this subsection we are
interested of the spectral properties of the operator
$\mathcal{L}$ defined by the
  \begin{equation}\label{3.12}
    \mathcal{L}=-\partial_x^2+c-\phi.
  \end{equation}
Let us first mention that (\ref{3.5}) reduces now to
$$X_0=c+\sqrt{c^2+2}, \quad
\Delta=\left(c-\sqrt{c^2+2}, c+2\sqrt{c^2+2}\right)
$$
By the definition of $a, b$ and $U(s)$ one obtains
$$\begin{array}{rl}
U(s)&\equiv \displaystyle \frac{1}{3}(\phi_0-s)[s^2+(\phi_0-3c)s
-(2\phi_0^2+3c\varphi_0-6\phi_2)]\\[4mm]
&\displaystyle=\frac{1}{3}(s-\phi_0)(\phi_1-s)(s+\phi_1+\phi_0-3c).
\end{array}$$
We note that the last equality is a consequence of Proposition 1,  which implies that $U(\phi_1)=U(\phi_0)=0$. To obtain an
explicit formula for the travelling wave $\phi_c$, we substitute
$\sigma=\phi_0+(\phi_1-\phi_0)z^2$, $z>0$ in order to express the
above integral  as an elliptic integral of the first kind in a
Legendre form. One obtains
$$\int_0^{Z(s)}\frac{dz}{\sqrt{(1-z^2)(\kappa'^2-k^2z^2)}}=\alpha(s-\xi),
$$
where
\begin{equation}\label{3.12a}
Z(s)=\sqrt{\frac{\phi_c(s)-\phi_0}{\phi_1-\phi_0}},\quad
k^2=\frac{\phi_1-\phi_0}{\phi_0+2\phi_1-3c}, \quad
\kappa^2+\kappa'^2=1, \quad
\alpha=\sqrt{\frac{\phi_0+2\phi_1-3c}{12}}.
\end{equation}
 Thus we
get the expression
\begin{equation}\label{3.13}
\phi_c(s)=\phi_0+(\phi_1-\phi_0)cn^2(\alpha(s-\xi);k).
\end{equation}
To calculate the period of $\phi_c$, we use (3.7) and the same
procedure as above. In this way we get
\begin{equation}
\label{a15}
T= 2\int_{\varphi_0}^{\varphi_1}\frac{d\sigma}{\sqrt{U(\sigma)}}=
 \frac{2}{\alpha}\int_0^1\frac{dz}{\sqrt{(1-z^2)(1-k^2z^2)}}=
 \frac{2K(k)}{\alpha}.
 \end{equation}

We return to the operator ${\mathcal L}$ defined by (\ref{3.12}),
where $\phi_c$ is determined by ({\ref{3.13}). Consider the
spectral problem
\begin{equation}\label{3.14}
\begin{array}{l}
\cL\psi=\lambda\psi,\\
\psi(0)=\psi(T),\;\psi'(0)=\psi'(T).\end{array}
\end{equation}
 We
will denote the operator just defined again by $\cL$. It is a
self-adjoint operator acting on $L^2_{per}[0,T]$ with
$D(\cL)=H^2([0,T])$. From the Floquet theory applied to
(\ref{3.14}) it follows} that its spectrum is purely discrete,
\begin{equation}\label{3.15}
\lambda_0<\lambda_1\leq\lambda_2<\lambda_3\leq\lambda_4<\ldots
\end{equation}
where $\lambda_0$ is always a simple eigenvalue. If $\psi_n(x)$ is
the eigenfunction corresponding to $\lambda_n$, then
\begin{equation}\label{3.16}
\begin{array}{l}
\psi_0 \;\;\mbox{\rm has no zeroes in}\;\; [0,T];\\
\psi_{2n+1},\;\psi_{2n+2} \;\;\mbox{\rm have each just} \;\; 2n+2
\;\;\mbox{\rm zeroes in}\;\; [0,T).
\end{array}
\end{equation}
\begin{proposition}
\label{prop:2}
 The linear
operator $\cL$ defined by $\eqref{3.14}$ has the following
spectral properties:
\begin{itemize}
\item[(i)]
  The first three eigenvalues of
$\cL$ are simple. \\
\item[(ii)] \it The second eigenvalue of $\cL$
is $\la_1=0$. \\
\item[(iii)] Remainder of the spectrum is consisted
by a discrete set eigenvalues.
\end{itemize}
\end{proposition}

\vspace{2ex} \noindent {\bf Proof.} By (\ref{3.2}), $\cL
\phi'_c=0$, hence  $\psi=\phi_c'$ is an eigenfunction
corresponding to zero eigenvalue. By Proposition 1 (iii) $\phi'$
has just two zeroes in $[0,T)$ and therefore by (\ref{3.16})
either $0=\lambda_1<\lambda_2$ or $\lambda_1<\lambda_2=0$ or
$\lambda_1=\lambda_2=0$. We are going to verify that only the
first possibility $0=\lambda_1<\lambda_2$ can occur. From the
definition of $k$ and $\alpha$ one obtains that
$$\phi_0+2\phi_1-3c=12\alpha^2,\quad
\phi_1-\phi_0=12k^2\alpha^2.$$ Then using ({\ref{3.13}) we get
$$\begin{array}{rl}
\cL & =-\partial_x^2+c-\phi_0-(\phi_1-\phi_0)cn^2(\alpha x;k)\\[2mm]
 & =-\partial_x^2+c-\phi_1+(\phi_1-\phi_0)sn^2(\alpha x;k)\\[2mm]
& = -\partial_x^2-\alpha^2[4k^2+4-12k^2sn^2(\alpha x;k)]\\[2mm]
& = \alpha^2[-\partial_y^2 -4k^2-4+12k^2sn^2(y;k)]\equiv
 \alpha^2\Lambda \end{array}$$
where $y=\alpha x$. The operator $\Lambda$ is related to Hill's
equation with Lam\'e potential
$$\Lambda w=-\frac{d^2}{dy^2}w+[12k^2sn^2(y;k)-4k^2-4]w=0$$
and its spectral properties in the interval $[0,2K(k)]$ are well
known \cite{ABS, HIK,In}. The first three (simple) eigenvalues and
corresponding periodic eigenfunctions of $\Lambda$ are
$$\begin{array}{l}
\mu_0=k^2-2-2\sqrt{1-k^2+4k^4}<0, \\
\psi_0(y)=dn(y;k)[1-(1+2k^2-\sqrt{1-k^2+4k^4})sn^2(y;k)]>0,\\[2mm]
\mu_1=0,\\
\psi_1(y)=dn(y;k)sn(y;k)cn(y;k)=\frac12(d/dy)sn^2(y;k),\\[2mm]
\mu_2=k^2-2+2\sqrt{1-k^2+4k^4}>0,\\
\psi_2(y)=dn(y;k)[1-(1+2k^2+\sqrt{1-k^2+4k^4})sn^2(y;k)].\end{array}$$
As the eigenvalues of $\cL$ and $\Lambda$ are related by
$\lambda_n=\alpha^2\mu_n$ we conclude that  the first three
eigenvalues of (\ref{3.14}) are simple and moreover $\lambda_0<0$,
$\lambda_1=0$, $\lambda_2>0$. The corresponding eigenfunctions are
$\psi_0(\alpha x)$, $\psi_1(\alpha x)=const.\phi'_c(x)$ and
$\psi_2(\alpha x)$. $\Box$

\subsection{The operator $\cL_{mKdV}$}
\label{sec:32}
Consider the modified Korteweg-de Vries equation
  \begin{equation}\label{mkdv}
    u_t+3u^2u_x+u_{xxx}=0.
   \end{equation}
   Traveling wave solutions in this case satisfy the equation
    \begin{equation}\label{3.21}
      -c\phi'+3\phi^2\phi'+\phi'''=0.
    \end{equation}
    Integrating yields
     \begin{equation}\label{3.22}
       \phi''=a+c\phi-\phi^3.
     \end{equation}
     We   consider  the "symmetric" case $a=0$ only.
     Integrating once again, we get
       \begin{equation}\label{3.23}
         \phi'^{2}=b+c\phi^2-{\frac{\phi^4}{2}}.
       \end{equation}

Hence the periodic solutions are given by the periodic
trajectories $H(\phi,\phi')=b$ of the Hamiltonian vector field
$dH=0$ where
$$H(x,y)=y^2+{\frac{x^4}{4}}-c{\frac{x^2}{2}}.$$

 Then
there are two possibilities

\vspace{1ex} \noindent 1.1) ({\it outer case}): for any $b>0$ the
orbit defined by $H(\phi,\phi')=b$ is periodic and oscillates
around the eight-shaped loop $H(\phi,\phi')=0$ through the saddle
at the origin.

\vspace{1ex} \noindent 1.2) ({\it left and right cases}): for any
$b\in(-\frac12c^2,0)$ there are two periodic orbits defined by
$H(\phi,\phi')=b$ (the left and right ones). These are located
inside the eight-shaped loop and oscillate around the centers at
$(\mp\sqrt{c},0)$, respectively.

 We will consider the left and right cases
of Duffing oscillator.

\vspace{2ex} In the left and the right cases, let us denote by
$\phi_1>\phi_0>0$ the positive roots of the quartic equation
${\frac{\phi^4}{2}}-a\phi^2-b=0$. Then, up to a translation, we
obtain the respective explicit formulas
\begin{equation}\label{3.24}
\phi(z)=\mp \phi_1 dn(\alpha z; k),\quad
k^2=\frac{\phi_1^2-\phi_0^2}{\phi_1^2}
=\frac{2\phi_1^2-2c}{\phi_1^2}, \quad
\alpha={\frac{\phi_1}{\sqrt{2}}}, \quad T=\frac{2K(k)}{\alpha}.
\end{equation}
Now
  \begin{equation}\label{3.24a}
  \mathcal{L}=-\partial_x^2+c-3\phi^2.
  \end{equation}
  We use (\ref{3.24})  to
rewrite the operator $\mathcal{L}$ in an  appropriate form. From the
expression for $\phi(x)$ from (\ref{3.24})
     and the relations
     between the elliptic functions $sn(x)$, $cn(x)$ and $dn(x)$, we obtain
       $$\mathcal{L}=\alpha^{2}[ -\partial_{y}^{2}+6k^{2} sn^{2}(y)-4-k^2] $$
     where $y=\alpha x$.

     It is well-known that the first five eigenvalues of
     $\Lambda =-\partial_{y}^{2}+6k^{2}sn^{2}(y, k)$,
     with periodic boundary conditions on $[0, 4K(k)]$, where
     $K(k)$ is the complete elliptic integral of the first kind, are
     simple. These eigenvalues, with their  corresponding eigenfunctions are as follows
      $$\begin{array}{ll}
         \nu_{0}=2+2k^2-2\sqrt{1-k^2+k^4},
         & \psi_{0}(y)=1-(1+k^2-\sqrt{1-k^{2}
         +k^{4}})sn^{2}(y, k),\\[1mm]
         \nu_{1}=1+k^{2}, & \psi_{1}(y)=cn(y, k)dn(y, k)
         =sn'(y, k),\\[1mm]
         \nu_{2}=1+4k^{2}, & \psi_{2}(y)=sn(y, k)dn(y, k)
         =-cn'(y, k),\\[1mm]
         \nu_{3}=4+k^{2}, & \psi_{3}(y)=sn(y, k)cn(y, k)
         =-k^{-2}dn'(y, k),\\[1mm]
         \nu_{4}=2+2k^{2}+2\sqrt{1-k^{2}+k^{4}},
         & \psi_{4}(y)=1-(1+k^{2}+\sqrt{1-k^{2}
         +k^{4}})sn^{2}(y, k).
        \end{array}
      $$
      It follows that the first three eigenvalues of the operator
      $\mathcal{L}$, equipped with periodic boundary condition on $[0,2K(k)]$
      (that is, in the case of left and right family),
      are simple and $\lambda_0=\alpha^2(\nu_0-\nu_3)<0, \;
      \lambda_1=\alpha^2(\nu_3-\nu_3)=0, \;
      \lambda_{2}=\alpha^2(\nu_4-\nu_3)>0$.
The corresponding eigenfunctions are $\chi_0=\psi_0(\alpha x),
\chi_1=\phi'(x), \chi_2=\psi_4(\alpha x)$.

Thus, we have proved the following
\begin{proposition}
\label{prop:3}
 The linear
operator $\cL$ defined by $\eqref{3.24a}$ has the following
spectral properties:
\begin{itemize}
\item[(i)]
  The first three eigenvalues of
$\cL$ are simple. \\
\item[(ii)] \it The second eigenvalue of $\cL$
is $\la_1=0$, which is simple. \\
\item[(iii)] The rest of the
spectrum   consists of a  discrete set of eigenvalues, which are strictly positive.
\end{itemize}
\end{proposition}

\section{Proof of Theorem \ref{theo:2}}
\label{sec:44}
We first consider the cases of the KdV and the modified KdV equations.
\subsection{The KdV and modified KdV equations}
We need to check the assumptions of Theorem \ref{theo:1} for the
operator $L(k)=-\p_x \cL\p_x+k^2$, where $\cL$ is either the
 operator associated to the KdV equation, constructed in Section \ref{sec:31} or the operator associated to the mKdV equation, constructed in Section \ref{sec:32}.

 Clearly,
$L(0)=-\p_x \cL\p_x$ is bounded from below self-adjoint operator,
so that its spectrum consists of eigenvalues with finite
multiplicities
$$
\si(L(0))=\la_0(L(0))\leq \la_1(L(0))\leq \ldots
$$
Thus, one may apply the Courant principle for the first
eigenvalue. We have
$$
\la_1(L(0))=\sup_{z\neq 0}\inf_{u\perp z} \f{\dpr{L(0)
u}{u}}{\|u\|^2}.
$$
and as a consequence the infimum in $u$ may be taken only on
functions with mean value zero.
  Taking $z=\psi'_0$ and the identity
  $\dpr{L(0)u}{u}=\dpr{-\p_x \cL \p_x u}{u}=\dpr{\cL u'}{u'}$ allows us to write
$$
\la_1(L(0))\geq \inf_{u\perp \psi'_0}
\f{\dpr{\cL u'}{u'}}{\|u\|^2}= \inf_{u'\perp \psi_0} \f{\dpr{\cL u'}{u'}}{\|u\|^2}
$$
since  $\dpr{u'}{\psi_0}=-\dpr{u}{\psi'_0}=0$. Now, observe that since in both $\cL_{KdV}$ and $\cL_{mKdV}$ we have that there is only a single and simple negative eigenvalue, it follows that
$\cL|_{\{\psi_0\}^\perp}\geq 0$, i.e. $\dpr{\cL v}{v} \geq 0$, whenever $v\perp \psi_0$.
In particular, if $u'\perp \psi_0$,
$$
\f{\dpr{\cL u'}{u'}}{\|u\|^2}\geq 0.
$$
Thus,
$$
\la_1(L(0))\geq  \inf_{u'\perp \psi_0} \f{\dpr{\cL u'}{u'}}{\|u\|^2}\geq 0.
$$
  Thus $\la_1(L(0))\geq 0$. On the other hand, we have
that $0$ is an eigenvalue for $L(0)$, because $L(0)\phi_c=-\p_x
\cL \p_x[ \phi_c]=-\p_x \cL \phi'_c=0$.

We will now show that  there is a negative eigenvalue for $L_{KdV}(0)$ and $L_{mKdV}(0)$.
We claim that this will be enough for the proof of Theorem \ref{theo:2}.

Indeed, if we succeed in showing  $\la_0(L(0))<0$, and since we
have established $0\in\si(L(0))$ and $\la_1(L(0))\geq 0$, it
follows that $\la_1(L(0))=0$. In particular $\la_0(L(0))$ is a
simple eigenvalue, hence verifying the first hypothesis of Theorem
\ref{theo:1} with $k_0^2:=-\la_0(L(0))$. Moreover, $L'(k_0)=2k_0
Id$ and hence, the second condition of Theorem \ref{theo:1} is
trivially satisfied as well.

Thus, it suffices to show
\begin{equation}
\label{a:1} \la_0(L(0))=\inf_{u: \|u\|=1}  \dpr{\cL u'}{u'} <0
\end{equation}

Here we present a sufficient condition for this to happen.
Namely, we construct  $u':=t_0 \psi_0-t_2 \psi_2$ for some
coefficient $t_0, t_2$ to be found momentarily. To that end, we
shall first need
\begin{equation}
\label{a:5} t_0 \int_0^T \psi_0(y) dy-t_2 \int_0^T \psi_2(y) dy=0
\end{equation}
to ensure that such a periodic function $u$ exists\footnote{in
which case, we simply define the non-trivially zero function
$u(x):=\int_0^x (t_0 \psi_0(y)-t_2 \psi_2(y)) dy $, which in view
of \eqref{a:5} is defined up to a multiplicative constant.}. Since
both $\int_0^T \psi_0(y) dy\neq 0, \int_0^T \psi_2(y) dy\neq 0$,
we conclude that we may select $t_0, t_2\neq 0$ and
\begin{equation}
\label{a100}
\f{t_0}{t_2}= \f{\int_0^T \psi_2(y) dy}{\int_0^T \psi_0(y) dy}.
\end{equation}
Next, we compute (using \eqref{a100})
\begin{eqnarray*}
\dpr{\cL u'}{u'} & = & \dpr{\cL (t_0\psi_0-t_2  \psi_2)}{(t_0
\psi_0-t_2 \psi_2)} =
t_0^2 \|\psi_0\|_{L^2}^2 \la_0(\cL)+t_2^2 \|\psi_2\|_{L^2}^2\la_2(\cL)= \\
&=&  t_2^2(\int_0^T \psi_2(y) dy)^2 \left( \f{\|\psi_0\|^2 \la_0(\cL)}{
(\int_0^T \psi_0(y) dy)^2}+ \f{\|\psi_2\|^2 \la_2(\cL)}{(\int_0^T \psi_2(y)
dy)^2} \right).
\end{eqnarray*}
Thus, it remains to check that under the conditions in Theorem
\ref{theo:1},  the following inequality holds true
\begin{equation}
\label{a:10}
\f{\|\psi_0\|^2 \la_0(\cL)}{(\int_0^T \psi_0(y) dy)^2}+
\f{\|\psi_2\|^2 \la_2(\cL)}{(\int_0^T \psi_2(y) dy)^2}<0.
\end{equation}
Thus, we have reduced the proof of Theorem \ref{theo:2} and Theorem \ref{theo:5} to checking \eqref{a:10} for $\cL_{KdV}$ and $\cL_{mKdV}$ respectively.

\subsubsection{Proof of \eqref{a:10} for $\cL_{KdV}$}

In the case of  Korteweg-de Vries equation using (\ref{3.12a}) and
  identities
    $$\begin{array}{ll}
      sn^2(x)={\frac{1}{\kappa^2}}(1-dn^2(x))\\
      \\
      \int_{0}^{K}{dn(x)}dx={\frac{\pi}{2}}\\
      \\
      \int_{0}^{K}{dn^3(x)}dx={\frac{\pi(2-\kappa^2)}{4}}\\
      \\
      \int_{0}^{K}{dn^2(x)sn^2(x)}dx={\frac{(2\kappa^2-1)E(\kappa)+(1-\kappa^2)K(\kappa)}{3\kappa^2}}\\
      \\
      \int_{0}^{K}{dn^2(x)sn^4(x)}dx={\frac{(8\kappa^4-3\kappa^2-2)E(\kappa)+2(1+\kappa^2-2\kappa^4)K(\kappa)}{15\kappa^4}}
      \end{array}$$
      we get\footnote{In the derivation of the formulas below,
      we have used the symbolic integration feature of the {\it Mathematica} software.}
      $$\begin{array}{ll}
        \int_{0}^{T}{\psi_0(\alpha x)}dx&={\frac{\pi}{\alpha}}\left[{\frac{2-\kappa^2}{2\kappa^2}}
        (1+2\kappa^2-\sqrt{1-\kappa^2+4\kappa^4})+{\frac{{\sqrt{1-\kappa^2+4\kappa^4}-1-\kappa^2}}{\kappa^2}}\right]\\
        \\
        \int_{0}^{T}{\psi_0^2(\alpha x)}dx&={\frac{2}{\alpha}}\left(
        E(\kappa)+{\frac{2(-1-2\kappa^2+\sqrt{1-\kappa^2+4\kappa^4})((-1+2\kappa^2)E(\kappa)-(-1+\kappa^2)K(\kappa))}{3\kappa^2}}\right.
        \\
        \\
        &\left.
        +{\frac{(-1-2\kappa^2+\sqrt{1-\kappa^2+4\kappa^4})^2((-2-3\kappa^2+8\kappa^4)E(\kappa)+2(1+\kappa^2-2\kappa^4)K(\kappa))}{15\kappa^4}}\right)\\
        \\
         \int_{0}^{T}{\psi_2(\alpha
         x)}dx&={\frac{\pi}{\alpha}}\left[{\frac{2-\kappa^2}{2\kappa^2}}(1+2\kappa^2+\sqrt{1-\kappa^2+4\kappa^4})
         -{\frac{1+\kappa^2+\sqrt{1-\kappa^2+4\kappa^4}}{\kappa^2}}\right]\\
         \\
         \int_{0}^{T}{\psi_2^2(\alpha x)}dx&={\frac{2}{\alpha}}\left(
        E(\kappa)+{\frac{2(1+2\kappa^2+\sqrt{1-\kappa^2+4\kappa^4})((1-2\kappa^2)E(\kappa)+(-1+\kappa^2)K(\kappa))}{3\kappa^2}}\right.
        \\
        \\
        &\left.
        +{\frac{(1+2\kappa^2+\sqrt{1-\kappa^2+4\kappa^4})^2((-2-3\kappa^2+8\kappa^4)E(\kappa)+2(1+\kappa^2-2\kappa^4)K(\kappa))}{15\kappa^4}}\right)
     \end{array}$$

     $$\begin{array}{ll}
        \lambda_0(\cL)=\alpha^2(k^2-2-2\sqrt{1-k^2+4k^4})\\
        \\
        \lambda_2(\cL)=\alpha^2(k^2-2+2\sqrt{1-k^2+4k^4}),
       \end{array}
      $$
      where $\alpha$ is given by (\ref{3.12a}).
      Thus, we have an explicit formula to work with in order to show \eqref{a:10}.
\begin{figure}[h]
\centering
\includegraphics[width=8cm,height=6cm]{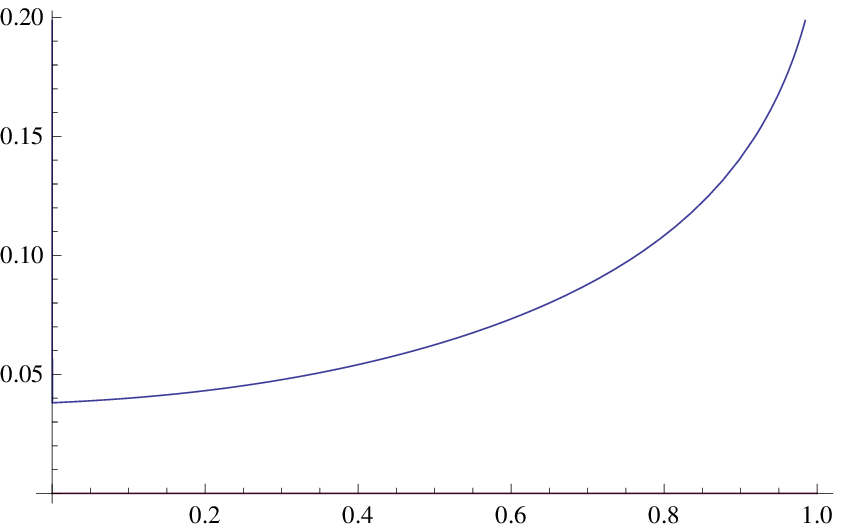}
\caption{This is a graph of the function
$h(\ka)= \f{\int_0^{2K(\ka)}\psi_2(x)dx}{\sqrt{\la_2} \|\psi_2\|}-
\f{\int_0^{2K(\ka)}\psi_0(x)dx}{\sqrt{|\la_0|} \|\psi_0\|} $.
Note that positivity of $h$ is equivalent to the validity of \eqref{a:10}.}
\label{fig1}
\end{figure}

From the graph in Figure \ref{fig1}, it is clear that the inequality \eqref{a:10} holds for all values of the parameter $\ka$.

\subsubsection{Proof of \eqref{a:10} for $\cL_{mKdV}$}
  In the case of Modified Korteweg-de Vries equation using (\ref{3.24}) and
  identities
    $$\begin{array}{ll}
      sn^2(x)={\frac{1}{\kappa^2}}(1-dn^2(x))\\
      \\
      \int_{0}^{K}{dn^2(x)}dx=E(\kappa)\\
      \\
      \int_{0}^{K}{sn^4(x)}dx={\frac{1}{3\kappa^4}}\left[(2+\kappa^2)K(\kappa)-2(1+\kappa^2)E(\kappa)\right]
      \end{array}$$

      we get that
      $$\begin{array}{ll}
        \int_{0}^{T}{\psi_0(\alpha x)}dx&={\frac{2}{\alpha \kappa^2}}(\sqrt{1-\kappa^2+\kappa^4}-1)K(\kappa)+(1+\kappa^2-\sqrt{1-\kappa^2+\kappa^4})E(\kappa)\\
        \\
        \int_{0}^{T}{\psi_0^2(\alpha x)}dx&={\frac{2}{\alpha}}\left( K(\kappa)-2(1+\kappa^2-\sqrt{1-\kappa^2+\kappa^4}){\frac{K(\kappa)-E(\kappa)}{\kappa^2}}\right. \\
        \\
        &\left. +(1+\kappa^2-\sqrt{1-\kappa^2+\kappa^4})^2
     {\frac{(2+\kappa^2)K(\kappa)-2(1+\kappa^2)E(\kappa)}{3\kappa^4}}\right) \\
         \\
         \int_{0}^{T}{\psi_2(\alpha
         x)}dx&={\frac{2}{\alpha \kappa^2}}(\sqrt{1-\kappa^2+\kappa^4}+1+\kappa^2)E(\kappa)-(1+\sqrt{1-\kappa^2+\kappa^4})K(\kappa)\\
         \\
         \int_{0}^{T}{\psi_2^2(\alpha x)}dx&={\frac{2}{\alpha}}\left( K(\kappa)-2(1+\kappa^2+\sqrt{1-\kappa^2+\kappa^4}){\frac{K(\kappa)_E(\kappa)}{\kappa^2}}\right. \\
         \\
         &\left. +(1+\kappa^2+\sqrt{1-\kappa^2+\kappa^4})^2
     {\frac{(2+\kappa^2)K(\kappa)-2(1+\kappa^2)E(\kappa)}{3\kappa^4}}\right)\\
      \\
       & \lambda_0=\alpha^2(k^2-2-2\sqrt{1-k^2+k^4})\\
        \\
        &\lambda_2=\alpha^2(k^{2}-2+2\sqrt{1-k^{2}+k^{4}}),
        \end{array}
      $$
      where $\alpha$ is given by (\ref{3.24}). Now
       the inequality \ref{a:10} is  is satisfies for all $\kappa
      \in (0,1)$. 
      Again, the graph below shows that the inequality \eqref{a:10} is satisfied for all values of the parameter $\ka$. 
      
\begin{figure}[h]
\centering
\includegraphics[width=8cm,height=6cm]{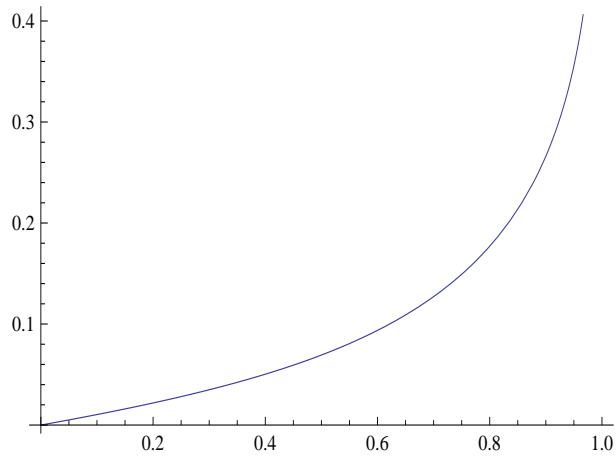}
\caption{This is a graph of the function
$h(\ka)= \f{\int_0^{2K(\ka)}\psi_2(x)dx}{\sqrt{\la_2} \|\psi_2\|}-
\f{\int_0^{2K(\ka)}\psi_0(x)dx}{\sqrt{|\la_0|} \|\psi_0\|} $.
Note that positivity of $h$ is equivalent to the validity of \eqref{a:10}.}
\label{fig2}
\end{figure}

\subsection{The nonlinear Schr\"odinger equation}
 In this section we will construct the periodic traveling wave
 solution for the quadratic and cubic nonlinear Schr\"odinger
 equations and investigate the spectral problems for corresponding
 operators. The results can be found in \cite{HIK3}, but for
 convenience we will present here.  
 We show that the matrix operator 
 $$
 \left(
 \begin{array}{cc}
 \cL_- & 0 \\
 0 & \cL_+
 \end{array}
 \right)
 $$ 
 has a single simple negative eigenvalue. The same will be true for the similar operator \\ $-J\left(
 \begin{array}{cc}
 \cL_- & 0 \\
 0 & \cL_+
 \end{array}
 \right) J=J^{-1} \left(
 \begin{array}{cc}
 \cL_- & 0 \\
 0 & \cL_+
 \end{array}
 \right) J $.  Thus, according to the instability criterium in Theorem 
 \ref{theo:1} and the representation \eqref{m3}, this implies that we can select a $k$ so that the operator $L(k)$ satisfies $(1)$ and $(2)$, whence we will have shown spectral instability. 
\subsubsection{Quadratic Schr\"odinger equation}

Consider the quadratic equation
\begin{equation}
\label{s1.1}
i u_t+u_{xx}+|u| u=0
\end{equation}
for a complex-valued function $u$.

  For $\varphi$ one obtains the equation \eqref{m1}, which is
    $$
    \varphi''-\omega \varphi +\varphi|\varphi|=0.
   $$
Therefore,
    \begin{equation}\label{s1.3}
    \varphi'^2-\omega \varphi^2 +\frac23\varphi^2|\varphi|=c
   \end{equation}
and $\varphi$ is periodic provided that the level set $H(x,y)=c$
of the Hamiltonian system $dH=0$, $$H(x,y)=y^2-\omega
x^2+\frac23x^2|x|,$$ contains a periodic trajectory (an oval). The
level set $H(x,y)=c$ contains two periodic trajectories if
$\omega>0$, $c\in(-\frac13\omega^3,0)$ and a unique periodic
trajectory if $\omega\in \mathbb{R}$, $c>0$. Under these
conditions, equation (\ref{s1.3}) becomes $H(\varphi,\varphi')=c$
and its solution $\varphi$ is periodic of period $T=T(\omega,c)$.

Below, we consider the case $c<0$. Then either $\varphi<0$ (the
left case) or $\varphi>0$ (the right case). To express $\varphi$
through elliptic functions, we denote by $\varphi_0>\varphi_1>0$
the positive solutions of $\frac23\rho^3-\omega \rho^2-c=0$. Then
$\varphi_1\leq |\varphi|\leq \varphi_0$ and one can rewrite
(\ref{s1.3}) as
\begin{equation}\label{s1.4}
\textstyle
\varphi'^2=\frac23(|\varphi|-\varphi_1)(\varphi_0-|\varphi|)(|\varphi|+\varphi_0+\varphi_1-\frac32\omega).
\end{equation}
Therefore
$2\varphi_0+\varphi_1>\varphi_0+2\varphi_1>\frac32\omega$.
Introducing a new variable $s\in(0,1)$ via
$|\varphi|=\varphi_1+(\varphi_0-\varphi_1)s^2$, we transform
(\ref{s1.4}) into
$$s'^2=\alpha^2(1-s^2)(k'^2+k^2s^2)$$ where $\alpha$, $k$, $k'$
are positive constants ($k^2+k'^2=1$) given by
$$\alpha^2=\frac{4\varphi_0+2\varphi_1-3\omega}{12}, \quad
k^2=\frac{2\varphi_0-2\varphi_1}{4\varphi_0+2\varphi_1-3\omega},
\quad
k'^2=\frac{2\varphi_0+4\varphi_1-3\omega}{4\varphi_0+2\varphi_1-3\omega}.$$
Therefore
\begin{equation}\label{s1.5}
|\varphi(x)|=\varphi_1+(\varphi_0-\varphi_1)cn^2(\alpha x;k).
\end{equation}

Consider in $[0,T]=[0, 2K(k)/\alpha]$ the differential
operators introduced earlier
\begin{equation}\label{s1.5a}
\cL_-=-\frac{d^2}{dx^2}+(\omega-2|\varphi|), \quad
\cL_+=-\frac{d^2}{dx^2}+(\omega-|\varphi|),
\end{equation}
 supplied with periodic
boundary conditions. By the above formulas,
$\varphi_0-\varphi_1=6\alpha^2k^2$,
 $2\varphi_0-\omega=4\alpha^2(1+k^2)$. Taking $y=\alpha x$ as an
independent variable in $\cL_-$, one obtains $\cL_-=\alpha^2\Lambda_1$
with an operator $\Lambda_1$  in $[0,2K(k)]$ given by
$$\begin{array}{rl} \Lambda_1&=\displaystyle -\frac{d^2}{dy^2}
+\alpha^{-2}[\omega-2(\varphi_1+(\varphi_0-\varphi_1)cn^2(y;k))]\\[5mm]
&\displaystyle= -\frac{d^2}{dy^2}+
\frac{\omega-2\varphi_0}{\alpha^2}+\frac{2(\varphi_0-\varphi_1)}{\alpha^2}sn^2(y;k)\\[5mm]
&\displaystyle= -\frac{d^2}{dy^2}-4(1+k^2)+12k^2sn^2(y;k).
\end{array}$$

The spectral properties of the operator  $\Lambda_1$  in
$[0,2K(k)]$ are well-known. The first three   eigenvalues are simple 
and moreover the  corresponding eigenfunctions of $\Lambda_1$ are given by 
  $$
  \begin{array}{ll}
    \mu_0=\kappa^2-2-2\sqrt{1-\kappa^2+4\kappa^4}<0 \\ 
    \psi_0(y)=dn(y;\kappa)[1-(1+2\kappa^2-\sqrt{1-\kappa^2+4\kappa^4})sn^2(y;\kappa)]>0\\ 
    \\
    \mu_1=0 \\ 
    \psi_1(y)=dn(y;\kappa)sn(y;\kappa)cn(y;\kappa)={\frac{1}{2}}{\frac{d}{dy}}cn^2(y;\kappa)\\
    \\
    \mu_2=\kappa^2-2+2\sqrt{1-\kappa^2+4\kappa^4}>0\\
    \psi_2(y)=dn(y;\kappa)[1-(1+2\kappa^2+\sqrt{1-\kappa^2+4\kappa^4})sn^2(y;\kappa)].
   \end{array}
  $$
  Since the eigenvalues of $\cL_-$ and $\Lambda_1$ are related via  
  $\lambda_n=\alpha^2 \mu_n$, it follows that the first three eigenvalues of the operator
      $\cL_-$, equipped with periodic boundary condition on $[0,2K(k)]$
      are simple and $\lambda_0<0, \lambda_1=0, \lambda_2>0$. The
      corresponding eigenfunctions are $\psi_0(\alpha x),
      \psi_1(\alpha x)=C \varphi'$ and $\psi_2(\alpha x)$.
In a similar way, since  $\cL_+=\alpha^2\Lambda_2$, one obtains that in 
$[0, 2K(k)]$
$$
\Lambda_2=-\frac{d^2}{dy^2}-2(1+k^2)+6k^2sn^2(y;k)+\omega/2\alpha^2.$$
To express $\omega$ through $\alpha$ and $k$, one should take into
account  the fact that in the cubic equation we used to determine
$\varphi_0$ and $\varphi_1$,
 we have that the coefficient at $\rho$  is zero. Therefore,
$$\textstyle \varphi_0\varphi_1+(\varphi_0+\varphi_1)(\frac32\omega-\varphi_0-\varphi_1)=0.$$ As
$\varphi_0=2\alpha^2+2\alpha^2k^2+\frac12\omega$,
$\varphi_1=2\alpha^2-4\alpha^2k^2+\frac12\omega$,  after replacing
these values in the above equation one obtains
$\omega^2=16\alpha^4(1-k^2+k^4)$. Since $\omega>0$, we finally
obtain
$$\Lambda_2=-\frac{d^2}{dy^2}+2(-1-k^2+\sqrt{1-k^2+k^4})+6k^2sn^2(y;k).$$
On the other hand, (\ref{s1.5}) yields
$$|\varphi|=2\alpha^2[1+k^2+\sqrt{1-k^2+k^4}-3k^2sn^2(y;k)].$$
 The first three eigenvalues and corresponding eigenfunctions
 of $\Lambda_2$ are as follows:
$$\begin{array}{ll} \lambda_0=0, & \psi_0=\varphi,\\[2mm]
 \lambda_1=2-k^2+2\sqrt{1-k^2+k^4}, & \psi_1=dn'(y;k)\\[2mm]
  \lambda_2=4\sqrt{1-k^2+k^4}, & \psi_2=1+k^2-\sqrt{1-k^2+k^4}-3k^2sn^2(y;k).
 \end{array}$$
The considerations above yield 
 \begin{proposition}
\label{prop:4}
 The linear
operator $\cL_-$ defined by $\eqref{s1.5a}$ has the following
spectral properties:
\begin{itemize}
\item[(i)]
  The first three eigenvalues of
$\cL_-$ are simple. \\
\item[(ii)] \it The second eigenvalue of $\cL_-$
is $\la_1=0$. \\
\item[(iii)] The rest of the spectrum of $\cL_-$ consists of 
 a discrete set of positive eigenvalues.
\end{itemize} 
The linear operator $\cL_+$ defined by $\eqref{s1.5a}$ has the
following spectral properties:
\begin{itemize}
\item[(i)]
 $\cL_+$  has no negative eigenvalue. \\
 \item[(ii)] \it The first eigenvalue of $\cL_+$
is zero, which is simple. \\ 
\item[(iii)]  The rest of the spectrum of $\cL_+$ consists of 
 a discrete set of positive eigenvalues.

\end{itemize}
\end{proposition}

 \subsubsection{Cubic Schr\"odinger equation}
 Consider the cubic nonlinear Schr\"odinger equation
   \begin{equation}\label{s2.1}
    iu_{t}+u_{xx}+|u|^{2}u=0,
   \end{equation}
  where $u=u(x,t)$ is a complex-valued function of $(x,t)\in \mathbb{R}^2$.

 For $\varphi$ one obtains the equation

    \begin{equation}\label{s2.2}
      \varphi''-\omega \varphi+\varphi^{3}=0.
    \end{equation}
    Integrating once again, we obtain
 \begin{equation}\label{s2.3}
 \varphi'^2-\omega \varphi^2+\frac12 \varphi^4=c
 \end{equation}
  and $\varphi$ is a periodic function
provided that the energy level set $H(x,y)=c$ of the Hamiltonian
system $dH=0$, $$H(x,y)=y^2-\omega x^2+\frac12x^4,$$ contains an
oval (a simple closed real curve free of critical points). The
level set $H(x,y)=c$ contains two periodic trajectories if
$\omega>0$, $c\in(-\frac12\omega^2,0)$ and a unique periodic
trajectory if $\omega\in \mathbb{R}$, $c>0$. Under these
conditions, the solution of (\ref{s2.2}) is determined by
$H(\varphi,\varphi')=c$ and $r$ is periodic of period
$T=T(\omega,c)$.

Below, we are going to consider the case $c<0$. Let us denote by
$\varphi_0>\varphi_1>0$ the positive roots of
${\frac{1}{2}}\varphi^4-\omega \varphi^2-c=0$. Then, up to a
translation, we obtain the respective explicit formulas
\begin{equation}\label{s2.4}
\varphi(z)=\mp \varphi_0 dn(\alpha z; k),\quad
k^2=\frac{\varphi_0^2-\varphi_1^2}{\varphi_0^2}
=\frac{-2\omega+2\varphi_0^2}{\varphi_0^2}, \quad
\alpha={\frac{\varphi_0}{\sqrt{2}}}, \quad T=\frac{2K(k)}{\alpha}.
\end{equation}
Here and below $K(k)$ and $E(k)$ are, as usual, the complete
elliptic integrals of the first and second kind in a Legendre
form. By (\ref{s2.4}), one also obtains $\omega=(2-k^2)\alpha^2$
and, finally,
\begin{equation}\label{s2.5}
T=\frac{2\sqrt{2-k^2}K(k)}{\sqrt{\omega}}, \quad k\in(0,1), \quad
T\in I=\left(\frac{2\pi}{\sqrt{\omega}},\infty\right).
\end{equation}

Again,  $\cL_-$ and $\cL_+$
       are given by
        \begin{equation}\label{s2.5a}
        \begin{array}{ll}
           \cL_-=-\partial_{x}^{2}+(\omega-3\varphi^{2}), \\[2mm]
           \cL_+=-\partial_{x}^{2}+(\omega -\varphi^{2}),
          \end{array}
         \end{equation}
         with periodic boundary conditions in $[0,T]$.

         \vspace{2ex}
     We use now (\ref{s2.4}) and (\ref{s2.5}) to rewrite operators
     $\cL_\pm$ in more
     appropriate form. From the expression for $\varphi(x)$ from (\ref{s2.4})
     and the relations
     between elliptic functions $sn(x)$, $cn(x)$ and $dn(x)$, we obtain
       $$\cL_-=\alpha^{2}[ -\partial_{y}^{2}+6k^{2} sn^{2}(y)-4-k^2] $$
     where $y=\alpha x$.

     It is well-known that the first five eigenvalues of
     $\Lambda_1=-\partial_{y}^{2}+6k^{2}sn^{2}(y, k)$,
     with periodic boundary conditions on $[0, 4K(k)]$, where
     $K(k)$ is the complete elliptic integral of the first kind, are
     simple. These eigenvalues and corresponding eigenfunctions are:
      $$\begin{array}{ll}
         \nu_{0}=2+2k^2-2\sqrt{1-k^2+k^4},
         & \phi_{0}(y)=1-(1+k^2-\sqrt{1-k^{2}
         +k^{4}})sn^{2}(y, k),\\[1mm]
         \nu_{1}=1+k^{2}, & \phi_{1}(y)=cn(y, k)dn(y, k)
         =sn'(y, k),\\[1mm]
         \nu_{2}=1+4k^{2}, & \phi_{2}(y)=sn(y, k)dn(y, k)
         =-cn'(y, k),\\[1mm]
         \nu_{3}=4+k^{2}, & \phi_{3}(y)=sn(y, k)cn(y, k)
         =-k^{-2}dn'(y, k),\\[1mm]
         \nu_{4}=2+2k^{2}+2\sqrt{1-k^{2}+k^{4}},
         & \phi_{4}(y)=1-(1+k^{2}+\sqrt{1-k^{2}
         +k^{4}})sn^{2}(y, k).
        \end{array}
      $$

      It follows that the first three eigenvalues of the operator
      $L_-$, equipped with periodic boundary condition on $[0,2K(k)]$
      (that is, in the case of left and right family),
      are simple and $\lambda_0=\alpha^2(\nu_0-\nu_3)<0, \;
      \lambda_1=\alpha^2(\nu_3-\nu_3)=0, \;
      \lambda_{2}=\alpha^2(\nu_4-\nu_3)>0$.
The corresponding eigenfunctions are $\psi_0=\phi_0(\alpha x),
\psi_1=\varphi'(x), \psi_2=\phi_4(\alpha x)$.

\vspace{2ex} Similarly, for the operator $\cL_+$ we have
$$\cL_+=\alpha^2[-\partial_y^2+2k^2sn^2(y, k)-k^2]$$
in the case of left and right family. The spectrum of
 $\Lambda_2=-\partial_y^2+2k^{2}sn^{2}(y, k)$ is formed
 by bands $[k^{2}, 1]\cup [1+k^{2}, +\infty)$. The
 first three eigenvalues and the corresponding eigenfunctions with
 periodic boundary conditions on $[0, 4K(k)]$ are simple and
       $$\begin{array}{ll}
          \epsilon_0=k^2, & \theta_0(y)=dn(y, k),\\[1mm]
          \epsilon_1=1, & \theta_1(y)=cn(y, k),\\[1mm]
          \epsilon_2=1+k^2, & \theta_2(y)=sn(y, k).
        \end{array}
      $$

From (\ref{s2.3}) it follows that zero is an eigenvalue of  $\cL_+$
and it is the first eigenvalue in the case of left and right
family, with corresponding eigenfunction $\varphi(x)$.

The above considerations gives an identical result to Proposition \ref{prop:4}  for the operators $\cL_\pm$ defined in \eqref{s2.5a}. Thus, in both the quadratic and cubic cases, we have obtained that there is a single negative eigenvalue for the matrix operator $\left(
\begin{array}{cc}
\cL_- & 0 \\
0 & \cL_+ 
\end{array}
\right)
$ 
and thus our proof is complete. 

\section{The defocusing modified KdV equation}
\label{ch:5}

Now consider the defocusing modified Korteweg-de Vries equation
  \begin{equation}\label{mkdv}
    u_t-3u^2u_x+u_{xxx}=0.
   \end{equation}
   We are looking for traveling wave solutions of the form $u(x,t)=\phi(x-ct), \; \; c<0$.

    So, if we substitute this specific solution in the defocusing mKdV and
consider the integration constant equal to zero then  $Q=Q_c$
satisfies the ordinary differential equation
\begin{equation}\label{dmkdv0}
\phi''-c\phi-\phi^3=0.
\end{equation}
 From this we obtain the first order differential equation
(in the associated quadrature form)
\begin{equation}\label{dmkdv1}
[\phi']^2=\frac{1}{2}(\phi^4+2c\phi^2+A),
\end{equation}
where $A$ is the integration constant and which need to be
different of zero for obtaining periodic profile solutions.
Analogously as in the case of modified Korteweg-de Vries equation,
we obtain the explicit form for the periodic traveling wave
solutions
$$
\phi_c(\xi)=\eta_2\operatorname{sn} (\alpha \xi;k),
$$
where $\eta_1>\eta_2>0$ are positive roots of the polynomial
$F(t)=t^4+2ct^2+A$ and $\alpha={\frac{\eta_1}{\sqrt{2}}}, \;
k^2=\eta_2^2/{\eta_1^2}\in (0,1)$. Since the function $sn(x)$ has
minimal period $4K(k)$ then the minimal period of $\phi$, $T$, is
given by $T=4K(k)/\alpha$. Moreover,
$$
k^2=\frac{-2c-\eta_1^2}{\eta_1^2}
$$

Now consider the spectral problem for the operator
$\mathcal{L}=-\partial_x^2+3\phi^2+c$,

\begin{proposition}\label{spdmkdv}
 Let $\phi$ be the snoidal wave solution of the defocusing Korteweg-de Vries equation. Let
\[
\lambda_0\leq\lambda_1\leq\lambda_2\leq\lambda_3\leq\lambda_4\leq
\cdot\cdot\cdot,
\]
connote the eigenvalues of the operator $\mathcal{L}$. Then
$$
\lambda_0< \lambda_1=0<\lambda_2<\lambda_3<\lambda_4
$$
are all simple whilst, for $j\geq 5$, the $\lambda_j$ are double
eigenvalues. The $\lambda_j$ only accumulate at $+\infty$.
\end{proposition}

\begin{proof}
 Since $\mathcal{L}\frac{d}{dx}\phi=0$ and
$\frac{d}{dx} \phi$ has 2 zeros in $[0,T)$, it follows that $0$ is
either $\lambda_1$ or $\lambda_2$. We will show that
$0=\lambda_1<\lambda_2$.

 From the expression for $\phi(x)$ , we obtain
       $$\mathcal{L}=\alpha^{2}[ -\partial_{y}^{2}+6k^{2} sn^{2}(y)-1-k^2] $$
     where $y=\alpha x$.

     The eigenvalues and corresponding eigenfunctions are:
      $$\begin{array}{ll}
         \nu_{0}=2+2k^2-2\sqrt{1-k^2+k^4},
         & \psi_{0}(y)=1-(1+k^2-\sqrt{1-k^{2}
         +k^{4}})sn^{2}(y, k),\\[1mm]
         \nu_{1}=1+k^{2}, & \psi_{1}(y)=cn(y, k)dn(y, k)
         =sn'(y, k),\\[1mm]
         \nu_{2}=1+4k^{2}, & \psi_{2}(y)=sn(y, k)dn(y, k)
         =-cn'(y, k),\\[1mm]
         \nu_{3}=4+k^{2}, & \psi_{3}(y)=sn(y, k)cn(y, k)
         =-k^{-2}dn'(y, k),\\[1mm]
         \nu_{4}=2+2k^{2}+2\sqrt{1-k^{2}+k^{4}},
         & \psi_{4}(y)=1-(1+k^{2}+\sqrt{1-k^{2}
         +k^{4}})sn^{2}(y, k).
        \end{array}
      $$

      It follows that the first five eigenvalues of the operator
      $\mathcal{L}$, equipped with periodic boundary condition on $[0,4K(k)]$
      are simple and $\lambda_0=\alpha^2(\nu_0-\nu_1)<0, \;
      \lambda_1=\alpha^2(\nu_1-\nu_1)=0, \;
      \lambda_{2}=\alpha^2(\nu_2-\nu_1)>0\; \lambda_3=\alpha^2(\nu_3-\nu_1)>0, \; \lambda_4=\alpha^2(\nu_4-\nu_1)>0$.
The corresponding eigenfunctions are $\chi_0=\psi_0(\alpha x),
\chi_1=\phi'(x), \chi_2=\psi_2(\alpha x), \; \chi_3=\psi_3(\alpha
x), \; \chi_4=\psi_4(\alpha x)$.

In the case of Defocusing Modified Korteweg-de Vries equation
 inequality \eqref{a:10} is
  equivalent to the inequality
 \begin{eqnarray*}
    & &  
    \frac{|
     (\sqrt{1-\kappa^2+\kappa^4}-1)K(\kappa)+(1+\kappa^2-\sqrt{1-\kappa^2+\kappa^4})E(\kappa)|}{\sqrt{|1+\kappa^2-2\sqrt{1-\kappa^2+\kappa^4}|}} < \\
     & < & 
 \frac{| (\sqrt{1-\kappa^2+\kappa^4}+1+\kappa^2)E(\kappa)-(1+\sqrt{1-\kappa^2+\kappa^4})K(\kappa)|}{\sqrt{1+\kappa^2+2\sqrt{1-\kappa^2+\kappa^4}}} 
  \end{eqnarray*}
 {\it However, one can see from the picture below, that this inequality does not hold for any value of $\ka$}. 
\begin{figure}[h]
\centering
\includegraphics[width=8cm,height=6cm]{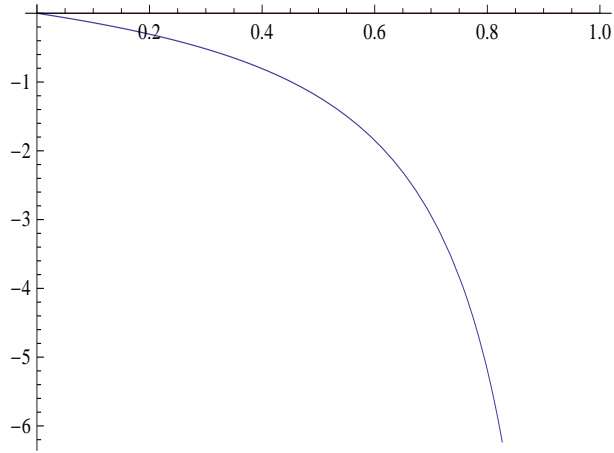}
\caption{Here, a plot of the difference of the two quantities is given. A positive function implies instability.}  
\label{fig3}
\end{figure}
Thus, our method fails to conclude the transversal instability of such waves. 
\end{proof}

\end{document}